\documentclass[journal,twoside,web]{ieeecolor}
\usepackage{tmi}

\usepackage{cite}
\usepackage{amsmath,amssymb,amsfonts}
\usepackage{algorithmic,  algorithm}
\usepackage{graphicx}
\usepackage{textcomp}
\usepackage{subfig}
\usepackage{bm, dsfont, url, braket}

\def\BibTeX{{\rm B\kern-.05em{\sc i\kern-.025em b}\kern-.08em
    T\kern-.1667em\lower.7ex\hbox{E}\kern-.125emX}}
\DeclareMathOperator*{\argmin}{arg\,min}

\title{Tuning-free multi-coil compressed sensing MRI with Parallel Variable Density Approximate Message Passing (P-VDAMP)}

\begin{document}
\author{Charles Millard, Mark Chiew, Jared Tanner, Aaron T Hess and Boris Mailh\'{e} 
\thanks{We are grateful to Thijs de Buck for sharing his angiogram data. This work was supported in part by an EPSRC Industrial CASE studentship with Siemens Healthineers under voucher 17000051,  in part by The Alan Turing
Institute under EPSRC grant EP/N510129/1, in part by EPSRC under grant EP/T013133/1, in part by the Royal Academy of Engineering under grant RF201617/16/23, and in part by the Wellcome Trust under grant 203139/Z/16/Z. The concepts and information presented in this paper are based on research results that are not commercially available.}
\thanks{Charles Millard, Mark Chiew and Aaron T. Hess are with the Wellcome Centre for Integrative Neuroimaging, FMRIB, Nuffield Department of Clinical Neurosciences, University of Oxford, Level 0, John Radcliffe Hospital, Oxford, OX3 9DU, UK. (email: charles.millard@ndcn.ox.ac.uk, mark.chiew@ndcn.ox.ac.uk and aaron.hess@ndcn.ox.ac.uk) }
\thanks{Jared Tanner is with the Mathematical Institute at the University of Oxford, Oxford, OX2 6GG, UK (email: tanner@maths.ox.ac.uk)}
\thanks{Boris Mailh\'{e} is with Siemens Healthineers, 755 College Rd E, Princeton,
NJ 08540, USA (email: boris.mailhe@siemens-healthineers.com)}
}

\maketitle 

\begin{abstract}
Magnetic Resonance Imaging (MRI) has excellent soft tissue contrast but is hindered by an inherently slow data acquisition process. Compressed sensing, which reconstructs sparse signals from incoherently sampled data, has been widely applied to accelerate MRI acquisitions. Compressed sensing MRI requires one or more model parameters to be tuned, which is usually done by hand, giving sub-optimal tuning in general. To address this issue, we build on previous work by the authors on the single-coil Variable Density Approximate Message Passing (VDAMP) algorithm, extending the framework to multiple receiver coils to propose the Parallel VDAMP  (P-VDAMP) algorithm. For Bernoulli random variable density sampling, P-VDAMP obeys a ``state evolution", where the intermediate per-iteration image estimate is distributed according to the ground truth corrupted by a zero-mean Gaussian vector with approximately known covariance. To our knowledge, P-VDAMP is the first algorithm for multi-coil MRI data that obeys a state evolution with accurately tracked parameters. We leverage state evolution to automatically tune sparse parameters on-the-fly with Stein's Unbiased Risk Estimate (SURE). P-VDAMP is evaluated on brain, knee and angiogram datasets and compared with four variants of the Fast Iterative Shrinkage-Thresholding algorithm (FISTA), including two tuning-free variants from the literature. The proposed method is found to have a similar reconstruction quality and time to convergence as FISTA with an optimally tuned sparse weighting and offers substantial robustness and reconstruction quality improvements over competing tuning-free methods. 
\end{abstract}

\begin{IEEEkeywords}
Magnetic Resonance Imaging (MRI), Compressed Sensing, Parallel Imaging, Approximate Message Passing, Stein's Unbiased Risk Estimate (SURE)
\end{IEEEkeywords}

\vspace{1cm}
\section{Introduction}

\IEEEPARstart{R}{educing}  acquisition time is one of the main challenges for rapid, dynamic Magnetic Resonance Imaging (MRI). Accordingly, there is considerable research attention on techniques that reconstruct sub-sampled k-space data without substantially sacrificing image quality. Parallel imaging  exploits the non-uniformity of coil sensitivites to reconstruct data that have been simultaneously acquired over a multichannel array of receiver coils \cite{ra1993fast, Pruessmann1999, Griswold2002, Uecker2014}. More recently, compressed sensing  \cite{Candes2006, Donoho2006} has been applied to MRI to accommodate further acceleration \cite{Lustig2007, Liang2009}, using presumed redundancy in the signal and regularized optimization techniques to reconstruct incoherently sampled data.

A frequently noted challenge for compressed sensing MRI (CS-MRI)  is how to tune model parameters \cite{Hsiao2012, Hollingsworth2014, Akasaka2016Optimization}, which is typically done by hand \cite{Jaspan2015}. The ideal parameter choice depends on a multitude of case-specific factors including the anatomy, underlying MR physics and sampling strategy, so hand-tuning generally yields sub-optimal reconstruction.  For reliably high quality CS-MRI reconstruction, a robust, accurate and computationally efficient method for automatic parameter tuning is desirable. 

Many existing methods for automatic parameter tuning depend on a testing procedure, where the reconstruction algorithm is run multiple times \cite{hansen1999curve, ramani2012regularization, giryes2011projected, weller2014monte}. In general this is a highly computationally intensive process, especially for models with multiple parameters. The focus of this paper is on algorithms that tunes model parameters on-the-fly, so that the algorithm only needs to be run once \cite{Khare2012, varela2021automatic}. 

A natural candidate for automatic on-the-fly parameter tuning is Stein's Unbiased Risk Estimate (SURE) \cite{Stein1981}. SURE efficiently estimates the squared error of a denoised vector, so can be used as a proxy for the actual error to tune parameters of a denoising function \cite{Donoho1995}. In CS-MRI, random sampling ensures that the aliasing is ``noise-like", suggesting that the application of SURE may be valid.  However, SURE only gives near-optimal parameter selection when the noise is drawn from a zero-mean Gaussian distribution, and the per-iteration aliasing of standard CS-MRI algorithms deviates considerably from a zero-mean Gaussian in practice. Despite no guarantee that the tuning is near-optimal, there has been empirical success employing SURE for tuning-free CS-MRI \cite{Ong2015, Khare2012}, suggesting that SURE is somewhat robust to deviations from the appropriate noise distribution.

This paper presents a solution to the problem of aliasing model mismatch via the Approximate Message Passing (AMP) framework\cite{Donoho2009, Ma2017, Rangan2019}. The AMP  algorithm was originally developed for sparse signals sampled with i.i.d. Gaussian sensing matrices\cite{Donoho2009}. For such sensing matrices, it is known that AMP obeys a ``state evolution", where the intermediate signal estimate is distributed according to the reference signal corrupted by white Gaussian aliasing. State evolution ensures that SURE's assumed aliasing distribution matches the actual aliasing, so automatic parameter tuning with SURE is truly near-optimal \cite{Mousavi2013, Guo2015}. However, when the sensing matrix does not satisfy its i.i.d. Gaussian assumption, such as for CS-MRI, AMP's state evolution may fail. The naive application of AMP to CS-MRI does not obey state evolution and does not perform well in practice. There has been a number of attempts to adapt AMP-based algorithms to CS-MRI including: Location Constrained AMP \cite{Sung2013}, Block-Matching 3D AMP-MRI \cite{Eksioglu2018, Metzler2016} Learned Denoising AMP \cite{9153368} and, for single coil measurements, Vector AMP (VAMP) for Image Reconstruction \cite{Rangan2019, Schniter2017a} and Damped Denoising VAMP \cite{Sarkar2020}. Although these modifications were found to empirically improve the convergence behaviour of AMP, they do not obey state evolution so cannot be leveraged for near-optimal parameter tuning with SURE. 

The authors' recently proposed single-coil VDAMP algorithm \cite{Milla2020, Millard2020a} is, to our knowledge, the first algorithm for Fourier sampled images that obeys a state evolution. Specifically, single-coil VDAMP obeys a ``colored state evolution", where the per-iteration aliasing  is accurately characterized by one real number per wavelet subband \cite{Milla2020}. This paper proposes the Parallel Variable Density Approximate Message Passing (P-VDAMP) algorithm, which extends single-coil VDAMP to multiple coils. For Bernoulli sampled k-space, where sampling is random and independent, P-VDAMP's intermediate aliasing is distributed according to a zero-mean Gaussian with a covariance matrix represented with a wavelet-domain vector with $N$ unique entries, where $N$ is the dimension of the reference image. We refer to this as as ``parallel colored state evolution". To our knowledge, P-VDAMP is the first algorithm for multi-coil MRI that obeys a state evolution with accurately tracked parameters, and therefore the first method where parameter tuning with SURE is truly near-optimal.


\section{Theory\label{sec:theory}}

Denote the $N$-dimensional data on coil $c$ as
\begin{align}
    \bm{y}_c &=  \bm{M}_\Omega (\bm{F}\bm{S}_c \bm{x}_0 + \bm{\varepsilon}_c), \hspace{1cm} c = \{1,2,\ldots,N_c\} \label{eqn:yc}
\end{align}
where $\bm{M}_\Omega$ is a diagonal undersampling mask with $1$ on the $j$th diagonal entry if $j \in \Omega$ and $0$ otherwise, where $\Omega$ is a sampling set with $|\Omega| = n$ for $n <N$.  The matrix $\bm{F}$ is the discrete Fourier transform and $\bm{S}_c$ is the diagonal sensitivity profile of coil $c$, which we assume have been normalized so that they satisfy $\sum_c \bm{S}_c^H\bm{S}_c = \mathds{1}_N$, and $\bm{x}_0$ is the image of interest.  The vector $\bm{\varepsilon}_c$ represents the measurement noise, which is modelled as complex Gaussian distributed with zero mean. 

The mathematical results in this section assume that the entries of the sampling set $\Omega$ are generated independently with arbitrary variable density $\mathrm{Pr}(j \in \Omega) = p_j$, referred to as ``Bernoulli sampling".   In other words, the sampling points correspond to Bernoulli distributed phase encodes with a fully sampled readout  dimension in a Cartesian sampling framework, and \eqref{eqn:yc} is a 2D subproblem that follows an inverse Fourier transform along the readout dimension. 

The experimental work in this paper can be reproduced with publicly available MATLAB code\footnote{\url{https://github.com/charlesmillard/P-VDAMP}}.

\subsection{Density compensation and error propagation \label{sec:err_prop}} 

In this section, we demonstrate that a Bernoulli sampled zero-filled image estimate has aliasing with a fully characterisable distribution. In other words, we show that state evolution holds at the first iteration. The mathematical tools developed in this section are preliminary to the full P-VDAMP algorithm presented in Section \ref{sec:statement_of_alg}. 

Consider the density compensated coil combined estimate from zero-filled data:
\begin{equation}
    \hat{\bm{x}}_0 =\sum_c \bm{S}_c^H \bm{F}^H \bm{P}^{-1}\bm{y}_c, \label{eqn:x_tilde}
\end{equation}
where $\bm{P} = \mathds{E}_{\Omega}\{ \bm{M}_\Omega \}$ is the expected value of the random sampling mask, so that $\bm{P}$ is diagonal with $j$th diagonal $p_j$. The $\bm{P}^{-1}$ matrix ensures that $\hat{\bm{x}}_0$ is an unbiased estimate of $\bm{x}_0$:
\begin{align*}
    \mathds{E}\{\hat{\bm{x}}_0\}     &=   \mathds{E}\{\sum_c \bm{S}_c^H \bm{F}^H \bm{P}^{-1}\bm{M}_\Omega (\bm{F}\bm{S}_c \bm{x}_0 + \bm{\varepsilon}_c)\} \\
    & = \sum_c \bm{S}_c^H \bm{F}^H \bm{P}^{-1} \bm{P} \bm{F}\bm{S}_c \bm{x}_0 = \bm{x}_0,
\end{align*}
where $ \mathds{E}\{\bm{M}_\Omega\} = \bm{P}$ and $ \mathds{E}\{\bm{\varepsilon}_c\} = \bm{0}$ has been used. 

It is often stated that the role of random sampling in CS-MRI is to ensure that the undersampling artifacts of the zero-filled estimate are ``noise-like". A primary result of this paper is to show that  when $\bm{M}_\Omega$ is Bernoulli random the undersampling artifacts of $\hat{\bm{x}}_0$ behave \textit{exactly} as Gaussian noise: $\hat{\bm{x}}_0= \bm{x}_0 + \mathcal{CN}(\bm{0}, \bm{\Sigma}^2_0)$, where $\mathcal{CN}(\bm{0}, \bm{\Sigma}^2)$ denotes the complex Gaussian distribution with zero mean and independent real and imaginary parts with covariance matrix $\bm{\Sigma}^2/2$. Moreover, we show that the covariance matrix $\bm{\Sigma}^2_0$ can be accurately and efficiently estimated with an $N$-dimensional vector, which we denote as $\bm{\tau}_0$. 

Image-domain aliasing is colored by the non-uniform spectral density of the image and the variable density sampling scheme \cite{Virtue2017}, so has a covariance matrix with strong off-diagonals. We employ an orthogonal wavelet transform to approximately decorrelate the aliasing, so that $\tilde{\bm{\Sigma}}^2_0$ is modelled as diagonal in the wavelet domain \cite{Milla2020}. Wavelets are commonly used in CS-MRI for sparsely representing an image, however, we emphasize that wavelets are used here as \textit{a tool for modelling the aliasing}, and that P-VDAMP is not constrained to wavelet-domain regularisation. The competitive performance of more sophisticated image models, including via neural networks, was demonstrated for single-coil VDAMP in \cite{Metzler2020}, although is not the focus of this paper. 

\begin{figure}
\centering
	  \includegraphics[width = 0.5\textwidth]{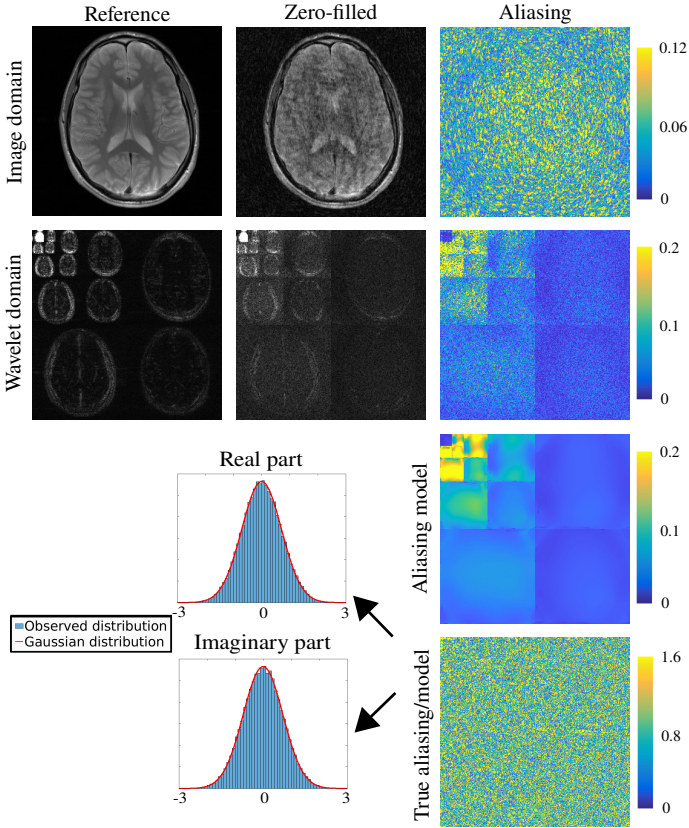}
      \caption{\label{fig:err_prop_example} The density compensated, zero-filled estimate $\hat{\bm{x}}_0$ for a $R = 10$ Bernoulli undersampled brain dataset. The image labelled ``Aliasing model" shows $\sqrt{\bm{\tau}_0}$ computed via \eqref{eqn:tau_imp}, and the image labelled ``True aliasing/model" shows $\bm{\Psi}(\bm{x}_0 - \hat{\bm{x}}_0)\oslash \sqrt{\bm{\tau}_0}$, where $\oslash$ is entry-wise division.  } 
\end{figure}

The appendix derives an approximation of the wavelet-domain covariance $\bm{\tau}_0 \approx \mathrm{diag}(\bm{\Psi} \tilde{\bm{\Sigma}}^2_0 \bm{\Psi}^H)$, where $ \bm{\Psi}$ is an orthogonal wavelet transform.  The $j$th entry is
\begin{equation}
	 \tau_{0,j} = \bm{\xi}_j^H \left( \sum_i |\hat{\Psi}_{ji}|^2\left[ \left(\frac{1 - p_i}{p_i}\right) \bm{y}_{i} \bm{y}_{i}^H + \bm{\Sigma}_{\varepsilon, i}^2 \right]\frac{m_i}{p_i} \right) {\bm{\xi}_{j}}. 
\label{eqn:tau_imp}
\end{equation}
Here, $\hat{\bm{\Psi}} = \bm{\Psi F}^H$,  $p_i$ and  $m_i$ are the $i$th diagonal of $\bm{P}$ and $\bm{M}_\Omega$ respectively, and $\bm{\Sigma}_{\varepsilon, i}^2 = \mathds{E} \{ \bm{\varepsilon}_i \bm{\varepsilon}_i^H \}$. The vectors $\bm{\xi}_j, \bm{y}_{i}$ and $\bm{\varepsilon}_i$ have dimensionality $N_c$; for instance, $\bm{y}_i$ is the vector of k-space location $i$ across all coils.  The notation $|\cdot|^2$ refers to the entry-wise magnitude squared.  As described in the appendix, the vector $\bm{\xi}_j$ approximates the coil-weighted spectral densities $\hat{\bm{\Psi}}^c = \bm{\Psi S}_c^H \bm{F}^H$ with  $\hat{\bm{\Psi}}_{ji}^c \approx \xi_{c,j} \hat{\bm{\Psi}}_{ji}$.  We suggest using 	$\xi_{c,j} = \braket{|\bm{\Psi}_j|^2, \mathrm{Diag}( \bm{S}_c)}$, where $\bm{\Psi}_j$ is the $j$th row of $\bm{\Psi}$. 


A visualisation of the error propagation formula is shown in Fig. \ref{fig:err_prop_example} for a $256\times 256$ brain dataset. Here, $N_c = 16$ coils sensitivities were estimated with ESPIRiT \cite{Uecker2014} via a $24\times 24$ central autocalibration region and compressed to $8$ virtual coils with Principle Component Analysis (PCA) \cite{Huang2008}, and an order 4 Daubechies wavelet $\bm{\Psi}$ at $4$ decomposition scales was employed. The histograms verify that the aliasing of $\hat{\bm{x}}_0$ is accurately approximated by $\bm{\tau}_0$, so that $\hat{\bm{x}}_0\approx \bm{x}_0 + \mathcal{CN}(\bm{0}, \bm{\Psi}\mathrm{Diag}(\bm{\tau}_0)\bm{\Psi}^H)$. 

The P-VDAMP algorithm described in Section \ref{sec:statement_of_alg} has such aliasing not only for the zero-filled initial estimate $\hat{\bm{x}}_0$ but also for all further iterations $\hat{\bm{x}}_k$: 
\begin{equation}
 \hat{\bm{x}}_k \approx \bm{x}_0 + \mathcal{CN}(\bm{0}, \bm{\Sigma}_k^2), \label{eqn:vdamp-se}
\end{equation}
where $\bm{\Sigma}_k^2 \approx \bm{\Psi}\mathrm{Diag}(\bm{\tau}_k)\bm{\Psi}^H$ and $\bm{\tau}_k$ is estimated with an error propagation formula based on \eqref{eqn:tau_imp}. We refer to   \eqref{eqn:vdamp-se} as P-VDAMP's ``parallel coloured state evolution".

\subsection{Automatic parameter tuning with SURE \label{sec:param_tuning_with_SURE}}

This section details how P-VDAMP's state evolution can be leveraged by applying SURE \cite{Stein1981}, building on previous work on tuning-free AMP in \cite{Mousavi2013, Guo2015, Bayati2013}. Consider a random vector $\bm{r} = \bm{w}_0 + \mathcal{CN}(\bm{0}, \mathrm{Diag}(\bm{\tau}))$, where $\bm{\tau}$ is known, and let an estimate of $\bm{w}_0$ be given by $\bm{f} (\bm{r}; \bm{\Theta})$, where $\bm{\Theta}$ is a parameter vector. Complex SURE takes the form \cite{Stein1981, Milla2020}
\begin{multline}
	cSURE(\bm{f}(\bm{r};\bm{\Theta}); \bm{\tau}) \\ =  \| \bm{f}(\bm{r}; \bm{\Theta}) - \bm{r} \|^2_2 + \bm{\tau}^T (2 \bm{\partial}(\bm{f}(\bm{r}; \bm{\Theta}))) - \bm{1}), \label{eqn:simpler_cSURE}
\end{multline}
where
\begin{equation}
   \partial_{i} ( \bm{f}(\bm{r}; \bm{\Theta})) =  \frac{1}{2}\left(\frac{\partial  \Re[f_i(\bm{r}; \bm{\Theta})]}{\partial \Re[r_i]}  + \frac{\partial \Im [f_i(\bm{r}; \bm{\Theta})]}{\partial \Im [r_i]} \right).
   \label{eqn:partial_def}
\end{equation}
Here, $\Re[\cdot]$ and $\Im[\cdot]$ are the real and imaginary parts respectively. cSURE is an unbiased estimate of the squared error of the denoised image, 
\begin{align}
	\mathds{E}\{\| \bm{f} (\bm{r}; \bm{\Theta}) - \bm{w}_0 \|^2_2 \} = \mathds{E}\{cSURE(\bm{f}(\bm{r};\bm{\Theta}); \bm{\tau}) \},
\end{align}
where $\mathds{E}\{ \cdot \}$ is the expectation over the random Gaussian noise. Therefore cSURE can be minimized as a proxy for minimizing the true squared error without requiring the ground truth $\bm{w}_0$, so a near-optimal denoiser satisfying
\begin{equation} 
	 \min_{\bm{\Theta}}  cSURE(\bm{f}(\bm{r};\bm{\Theta}); \bm{\tau})  \approx \min_{\bm{\Theta}} \| \bm{f}(\bm{r}; \bm{\Theta}) - \bm{w}_0 \|^2_2 \label{eqn:SUREprox}
\end{equation}
can be employed. 

In this paper, as in \cite{Milla2020}, we use the aliasing model \eqref{eqn:tau_imp} in conjunction with cSURE to employ soft thresholding with a per-subband threshold. To tune such a model by hand would be highly impractical; for instance, per-subband thresholding of a wavelet transform with $4$ decomposition scales would require 13 parameters to be tuned by hand.


\subsection{Statement of algorithm \label{sec:statement_of_alg}}

\begin{algorithm}[t]
\caption{P-VDAMP }
\textbf{Require:} Sampling set $\Omega$, coil sensitivities $\bm{S}_c$, wavelet transform $\bm{\Psi}$, probability matrix $\bm{P}$, measurements $\bm{y}_c$, measurement noise covariance $\bm{\Sigma}_{\varepsilon}^2$, denoiser $\bm{g}(\bm{r}; \bm{\tau})$, number of iterations $K_{it}$.
\begin{algorithmic}[1]
\STATE Compute $|\hat{\bm{\Psi}}|^2 = |\bm{\Psi}\bm{F}^H|^2$ and $\xi_{c,j}= \braket{|\bm{\Psi}_j|^2, \mathrm{Diag}( \bm{S}_c^H)}$ and set $\widetilde{\bm{r}}_0 = \bm{0}_N$
\FOR {$k =0,1,\ldots, K_{it}-1$} 
\STATE $\bm{z}_{k,c} = \bm{y}_c - \bm{M}_\Omega \bm{F S}_c \bm{\Psi}^H \widetilde{\bm{r}}_k$ \label{algline:3}
\STATE $\bm{r}_k = \widetilde{\bm{r}}_k + \bm{\Psi} \sum_c \bm{S}_c^H \bm{F}^H \bm{P}^{-1}\bm{z}_{k, c}$ \label{algline:4}
\STATE $  \tau_{k, j} =   \bm{\xi}_{j}^H \left( \sum_i |\hat{\Psi}_{ji}|^2 \left[ \left(\frac{1 - p_i}{p_i}\right) \bm{z}_{k, i} \bm{z}_{k, i}^H + \bm{\Sigma}_{\varepsilon, i}^2 \right]\frac{m_i}{p_i} \right) \bm{\xi}_{j}  $
\STATE $\hat{\bm{w}}_{k} = \bm{g}(\bm{r}_{k}; \bm{\tau}_{k})$ \label{algline:6}
\STATE $\bm{\alpha}_{k} = \braket{\bm{\partial} (\bm{g}(\bm{r}_{k}; \bm{\tau}_{k}))}_\mathrm{sband}$ \label{algline:7}
\STATE $\bm{c}_k = \bm{1}_N \oslash (\bm{1}_N - \bm{\alpha}_{k})$ \label{algline:8}
\STATE $\widetilde{\bm{r}}_{k+1} = \bm{c}_k \odot (\hat{\bm{w}}_{k} - \bm{\alpha}_k \odot \bm{r}_{k})$ \label{algline:9}
\ENDFOR
\RETURN $\hat{\bm{x}} = \bm{\Psi}^H\hat{\bm{w}}_k +\sum_c \bm{S}_c^H \bm{F}^H (\bm{y}_c - \bm{M}_\Omega \bm{F S}_c \bm{\Psi}^H \hat{\bm{w}}_k )$ \\ \hspace{0.64cm} or $\hat{\bm{x}} = \bm{\Psi}^H \bm{r}_k$
\end{algorithmic}
\label{alg:mult_coil_VDAMP}
\end{algorithm} 
 

The P-VDAMP algorithm is shown in Algorithm \ref{alg:mult_coil_VDAMP}. Here, vectors with subscript $c$ ($\bm{y}_c$ and $\bm{z}_{k,c}$) refer to an $N$-dimensional vector on the $c$th coil, and vectors with subscript $i$ or $j$ ($\bm{\xi}_{j}$ and $\bm{z}_{k,i}$) refer to the $N_c$-dimensional vector at the $i$th or $j$th k-space coefficient. In lines 3 and 5, the computation is for all coils $c$ and k-space locations $j$ respectively.

P-VDAMP is closely related to the single-coil VDAMP algorithm \cite{Milla2020}. The crucial difference between VDAMP and P-VDAMP is the aliasing model update in line 5. While VDAMP represents the aliasing with one real number per wavelet subband, P-VDAMP uses \eqref{eqn:tau_imp} to represent the aliasing with an $N$-dimensional vector, modelling the intra-subband modulation induced by the coil sensitivity profiles.

Lines 3-4,  referred to as ``density compensated gradient descent", is the multi-coil version of lines 3-4 of VDAMP \cite{Milla2020}. This can be understood as gradient descent with a modified forward model, so that \eqref{eqn:yc} is left-multiplied by $\bm{P}^{-1/2}$,
\begin{align}
    \bm{P}^{-1/2} \bm{y}_c &=  \bm{P}^{-1/2} \bm{M}_\Omega (\bm{F}\bm{S}_c \bm{x}_0 + \bm{\varepsilon}_c),
\end{align}
for all coils. This weighting ensures that the columns of the multi-coil sensing matrix are normalized in expectation over the random mask, a necessary component of state evolution. At the first iteration $k=0$, where $\widetilde{\bm{r}}_0 = \bm{0}$, lines 3-4 correspond to the density-compensated zero-filled estimate stated in  \eqref{eqn:x_tilde}. Using density compensated gradient descent implies that the data consistency term in P-VDAMP's cost function is a weighted $\ell_2$ norm, $\sum_c\|\bm{P}^{-1/2}(\bm{y}_c - \bm{M}_\Omega \bm{F}\bm{S}_c \bm{x})\|_2^2$, the consequences of which are discussed in Section \ref{sec:recon_qual}.

Line 5 applies the aliasing model update stated in \eqref{eqn:tau_imp} to the vector $\bm{z}_k$, which estimates the wavelet-domain covariance matrix of the aliasing of $\bm{r}_k$. Lines 6-9 are identical to lines 6-9 of the single-coil VDAMP algorithm \cite{Milla2020}. In line 6, a function $\bm{g}(\bm{r}_k; \bm{\tau}_k)$ is applied to the wavelet-domain estimate $\bm{r}_k$. Since, by state evolution, $\bm{r}_k$ is the Gaussian corrupted ground truth, $\bm{g}(\bm{r}_k; \bm{\tau}_k)$ is referred to as a \textit{denoiser}. In this work,  the denoiser was soft thresholding with per-subband soft thresholds tuned with  cSURE via the aliasing model $\bm{\tau}_k$, as described in Section \ref{sec:param_tuning_with_SURE}. Lines 7-9 correspond to a correction of the denoiser referred to in the AMP literature as ``Onsager correction", which removes the correlation that causes state evolution to fail for standard CS-MRI algorithms \cite{Donoho2009, Ma2017}. In line 7, the notation $\braket{\cdot}_{\mathrm{sband}}$ averages subbands, so that $\bm{\alpha}_{k}$ has the structure 
\begin{equation}
 \bm{\alpha}_k = \begin{bmatrix}
    \alpha_{k,1} \bm{1}_{N_1} \\
    \alpha_{k,2} \bm{1}_{N_2}\\
    \vdots \\
    \alpha_{k,1+3s} \bm{1}_{N_{1+3s}}
    \end{bmatrix},  \label{eqn:vecalpha}
\end{equation}
with
\begin{equation*}
    \alpha_{k,b} = \frac{1}{N_b} \sum_{j \in J_b}  \partial_j (\bm{g}(\bm{r}_{k}; \bm{\tau}_{k})),
\end{equation*}
where $J_b$ is the set of indices associated with subband $b$, $N_b = |J_b|$ and $\bm{\partial}(\cdot)$ is defined in \eqref{eqn:partial_def}. The notation $\oslash$ in line 8 refers to entry-wise division and $\odot$ in line 9 refers to entry-wise multiplication.


Line 11 of P-VDAMP gives two options for the algorithm output. The first uses gradient descent without density compensation.  We also suggest $\bm{\Psi}^H\bm{r}_k$, which is unbiased when state evolution holds. The qualitative and quanitative differences between these outputs are discussed in Section \ref{sec:num_exp}.

\subsection{Comments on P-VDAMP in practice \label{sec:mc_in_practice}}

We have found that some damping aids P-VDAMP's convergence, so that the denoiser is damped by a constant factor $0<\rho \leq 1$, as in \cite{Rangan2014, Rangan2019}. For $k >0$, we suggest replacing lines 6 and 7 of P-VDAMP, Algorithm \ref{alg:mult_coil_VDAMP}, with
\begin{align*}
	\hat{\bm{w}}_k & = \rho \bm{g}(\bm{r}_k; \bm{\tau}_k) + (1-\rho) \bm{\hat{w}}_{k-1}  \\
	\bm{\alpha}_k & = \rho \braket{\bm{\partial} (\bm{g}(\bm{r}_k; \bm{\tau}_k))}_{\text{sband}}.
\end{align*} 
Note that $\rho = 1$ corresponds to zero damping. In all experiments in this paper, the damping factor was set to $\rho = 0.75$. 

Rather than using a fixed number of iterations, we used a more practical stopping criterion based on the aliasing model $\bm{\tau}_k$. \cite{Metzler2020} In the experiments in this paper, P-VDAMP stops iterating when the mean-squared of the aliasing is predicted to have increased, $\braket{\bm{\tau}_k}>\braket{\bm{\tau}_{k-1}}$, or when $\braket{\bm{\tau}_k}$ changes by a small amount: $|\braket{\bm{\tau}_k} - \braket{\bm{\tau}_{k-1}}|/\braket{\bm{\tau}_{k-1}} < \epsilon$. We used  $\epsilon = 10^{-3}$ throughout this paper.

\section{Methods}

\subsection{Description of data}
P-VDAMP was evaluated on three types of MRI data: 
\begin{itemize}
	\item \textbf{Brain.} $256 \times 256$ acquisitions acquired on $N_c = 16$ coils, from an on-site Siemens 3T scanner. 
	\item \textbf{Knee.} $368 \times 640$ acquisitions acquired on $N_c = 15$ coils, from the publicly available fastMRI dataset \cite{zbontar2018fastmri}, which includes a mixture of field strengths. 
	\item \textbf{Time-of-flight Angiography.} $640 \times 506 \times 56$ acquisitions acquired on $N_c = 32$ coils, from an on-site Siemens 7T scanner, decomposed in to a series of 2D problems by taking an inverse Fourier transform along the first dimension. 
\end{itemize}
For each of the three categories, we used two fully sampled datasets. The two knee datasets were randomly selected from the fastMRI dataset. 

\subsection{Experimental method}

\begin{figure}
\centering
\subfloat[Sampling density $\bm{P}$ ]{%
    \includegraphics[width=0.4\columnwidth]{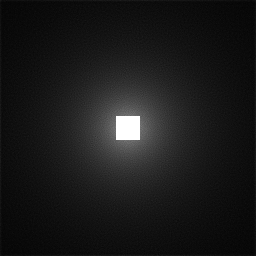}}
\hspace{0.3cm}
\subfloat[Random mask $\bm{M}_\Omega$]{%
    \includegraphics[width=0.4 \columnwidth]{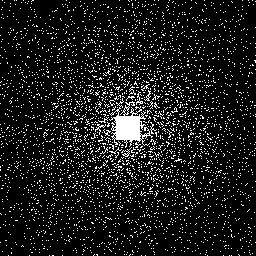}}
\caption{The sampling density and Bernoulli random mask employed in this paper for a $256 \times 256$ k-space at $R=10$. \label{fig:mask_examples}}
\end{figure}

All datasets were retrospectively undersampled at factors $R = 5, 10$. with variable density Bernoulli masks. The masks had $p_j =1$ in a $24 \times 24$ central autocalibration region used to estimate the coil sensitivities $\bm{S}_c$ with ESPIRiT \cite{Uecker2014} via the BART toolbox \cite{uecker2016bart}.  An example of the sampling probability and instance of a corresponding mask for a $256 \times 256$ k-space is shown in Fig. \ref{fig:mask_examples}.  

To reduce the reconstruction time, the coils were compressed to $8$ virtual coils with PCA \cite{Huang2008}: equation \eqref{eqn:tau_imp} is linear in $N$ but quadratic in $N_c$, so P-VDAMP especially benefits from coil compression. For all algorithms we used the order 4 Daubechies mother wavelet at $4$ decompositions scales for $\bm{\Psi}$. P-VDAMP's sensitivity to measurement noise is evaluated in Section \ref{sec:meas_noise_effect}; otherwise, the covariance $\bm{\Sigma}_{\varepsilon}^2$ was set to a matrix of zeros.

\subsection{Comparative algorithms}

We evaluated the reconstruction quality of both output options in line 11 of P-VDAMP, Algorithm \ref{alg:mult_coil_VDAMP}, referring to the $\bm{\Psi}^H\bm{r}_k$ output as ``Unbiased P-VDAMP", and the alternative as simply ``P-VDAMP". 

For comparative hand-tuned methods,  we implemented the Fast Iterative Shrinkage-Thresholding algorithm (FISTA) \cite{Beck2009} with a  $\lambda$  tuned using the reference to be approximately MSE optimal, referred to as ``Optimal FISTA", using an exhaustive search within a reasonable range.  Secondly, to reflect the more realistic case where the sparse parameter is not optimally tuned, we also ran FISTA with a sparse weighting shared across all reconstructions, which we refer to as ``Shared FISTA". For the shared sparse weighting we used the mean of the optimally tuned sparse weightings.

For comparative tuning-free methods, we ran FISTA with a per-subband threshold tuned automatically with SURE using a white alising model, referred to as ``SURE-IT"\cite{Khare2012}. We also evaluated a recently proposed parameter-free variant of FISTA, which automatically selects a per-scale sparse weighting based on k-means clustering of the zero-filled estimate, referred to in this work as ``Automatic-FISTA" (A-FISTA)\cite{varela2021automatic}. 

The stopping criteria employed for Optimal FISTA, Shared FISTA and SURE-IT was analogous to the P-VDAMP criteria described Section \ref{sec:mc_in_practice}. For A-FISTA, we found that such a stopping criterion typically stopped the algorithm prematurely, so to give A-FISTA the best chance of a high-quality reconstruction we used a fixed, large number of iterations. As in  \cite{varela2021automatic}, we chose 200 iterations. 

\begin{figure*}[th]
	\includegraphics[width = 2\columnwidth]{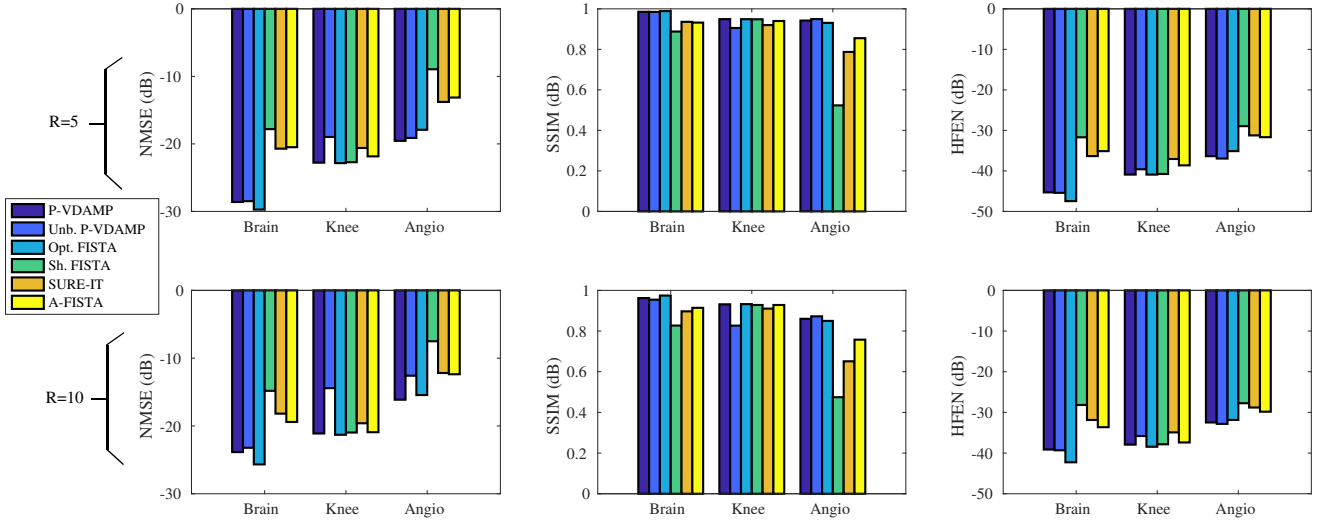}
      \caption{Reconstruction quality for all algorithms, where  the result for each data type has been averaged over the two datasets. \label{fig:mc_results}}
\end{figure*}

\subsection{Performance metrics \label{sec:perform_metric}}

Prior to computing performance metrics, we masked the images to mitigate background effects. For the brain and knee datasets, all pixels with absolute value less than 5\% of the maximum of the reference were masked. We also skull-stripped the angiograms with a hand-selected mask \cite{kalavathi2016methods} and computed the maximum projection along the readout dimension.

To quantify reconstruction quality  we used the Normalized Mean-Squared Error (NMSE), Structural Similarity Index Metric (SSIM)  \cite{Wang2004} and High-Frequency Error Norm (HFEN) \cite{ravishankar2010mr}, which quantifies the reconstruction quality of the image's fine features and is computed by applying a $15 \times 15$ Laplacian of Gaussian filter with standard deviation $1.5$ to $\hat{\bm{x}} - \bm{x}_0$ and computing the $\ell_2$ norm. We normalized the HFEN by dividing by the $\ell_2$ norm of the reference image. For the angiograms, we computed the SSIM after applying a vessel mask  computed by thresholding with Otsu's method \cite{otsu1979threshold, blob}.   %

\begin{figure*}
	\includegraphics[width=2\columnwidth]{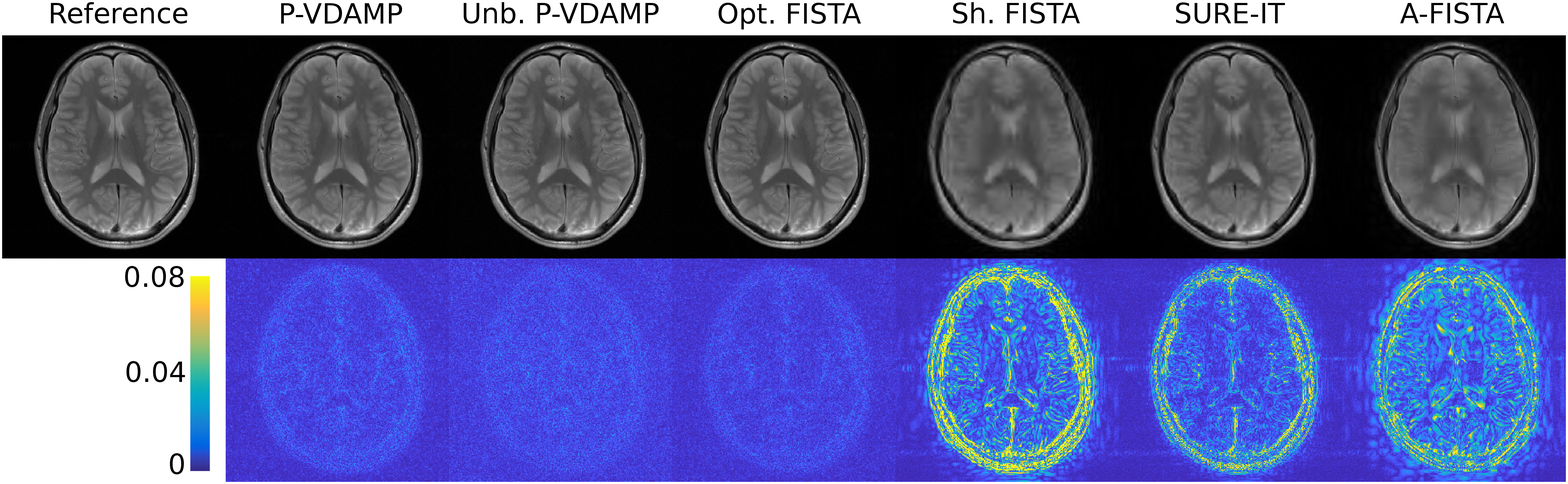}
	\caption{Reconstructions of a brain dataset sampled at $R = 5$ with error maps. Shared FISTA's sparse weighting is around a factor of 10 larger than for Optimal FISTA, which led to significant reconstruction artifacts. P-VDAMP performs competitively with Optimal FISTA and substantially better than the competing parameter-free methods SURE-IT and A-FISTA. The scale on the colorbar indicates that its maximum is 0.08 times the maximum of the reference image. \label{fig:brain_recon}}
	\vspace{0.5cm}
	\includegraphics[width=2\columnwidth]{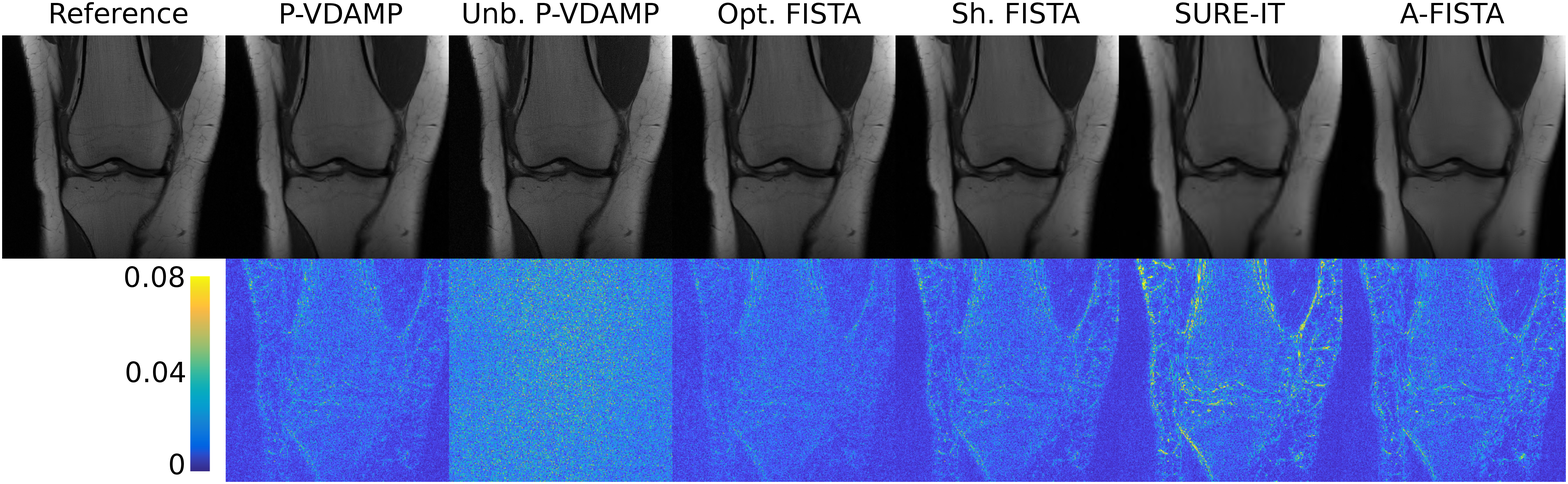}
	\caption{Reconstructions of a knee dataset sampled at $R = 10$, with error maps, cropped to a central $320 \times 320$ region. Unbiased P-VDAMP's NMSE is 2.5dB higher than the next worst, SURE-IT, but has competitive HFEN and performs well qualitatively. The competing parameter-free algorithms SURE-IT and A-FISTA smooth over details fine details in the bone. \label{fig:knee_recon}}
	\vspace{0.5cm}
	\includegraphics[width=2\columnwidth]{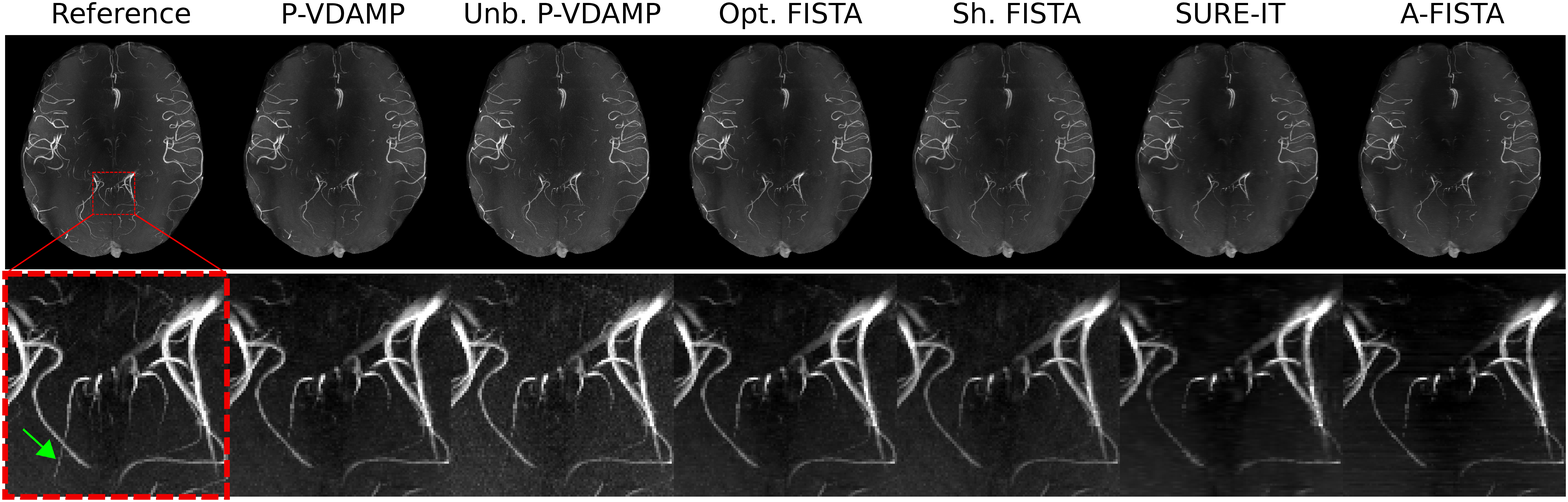}
	\caption{An angiogram sampled at $R = 5$, with a magnified $100 \times 100$ region.  The green arrow shows a challenging narrow blood vessel that is retained well for Unbiased P-VDAMP but mostly lost in competing methods. \label{fig:angio_recon}}
\end{figure*}

\begin{figure}
\centering
	\includegraphics[width =  \columnwidth]{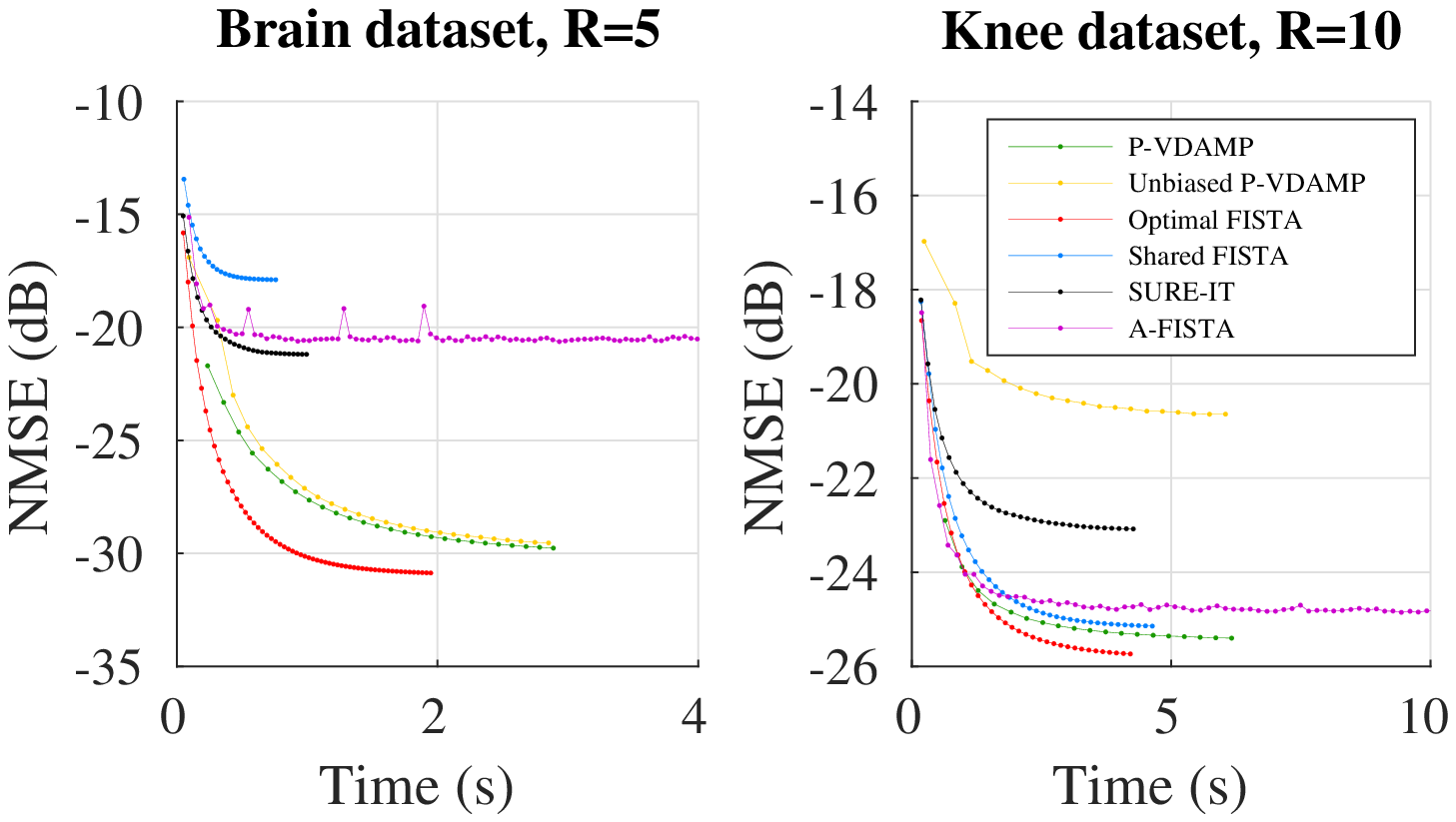}
	\caption{NMSE vs time for a brain and knee dataset at $R = 5$ and $R = 10$ respectively.  \label{fig:NMSE_vs_it_mc}}
\end{figure}

\section{Results and Discussion \label{sec:num_exp} }
\subsection{Reconstruction quality \label{sec:recon_qual}}

The NMSE, HFEN and SSIM of each is algorithm is shown  in Fig. \ref{fig:mc_results}. P-VDAMP consistently outperformed the competing parameter-free methods SURE-IT and A-FISTA. The improvement over SURE-IT emphasizes that P-VDAMP's state evolution is key for high-quality  threshold selection with SURE.  A-FISTA's performace is quite varied, which indicates that the quality of its parameter tuning is not as robust as P-VDAMP. The highly mixed performance of Shared FISTA shows that that the choice of sparse weighting is key, and that sharing its value for different types of data cannot be expected to yield a consistently high quality reconstruction. P-VDAMP peformed competitively with Optimal FISTA for all reconstructions, indicating that it is robust to the forward model and anatomy.  Since the fully sampled reference is not available for parameter tuning in prospectively undersampled scans, FISTA's parameter tuning cannot be expected to be near-optimal, so P-VDAMP is likely to outperform FISTA in practice.   Unbiased P-VDAMP's NMSE and SSIM was generally poorer than P-VDAMP, but had comparable HFEN.

Since P-VDAMP's sparse model has one parameter per subband while Optimal FISTA has a single, global parameter, one might expect that P-VDAMP would consistently outperform Optimal FISTA. There are two reasons why P-VDAMP does not always outperform Optimal FISTA. Firstly, P-VDAMP's $\bm{P}^{-1}$ in lines 3-4 of Algorithm \ref{alg:mult_coil_VDAMP},  which is essential for state evolution, is not advantageous from a reconstruction error perpective as it boosts errors at high frequencies, where the sampling probability is low. P-VDAMP outperforms Optimal FISTA when the benefits of its richer sparse model ``wins" over the disadvantage of needing to include $\bm{P}^{-1}$. The second reason is that damping causes P-VDAMP's aliasing model $\bm{\tau}_k$ to be imperfect, especially for high $k$, as discussed in Section \ref{sec:mc_approx_se}, so parameter selection with SURE is not always near-optimal in practice. 

\begin{figure}[t]
\centering
\includegraphics[width=1 \columnwidth]{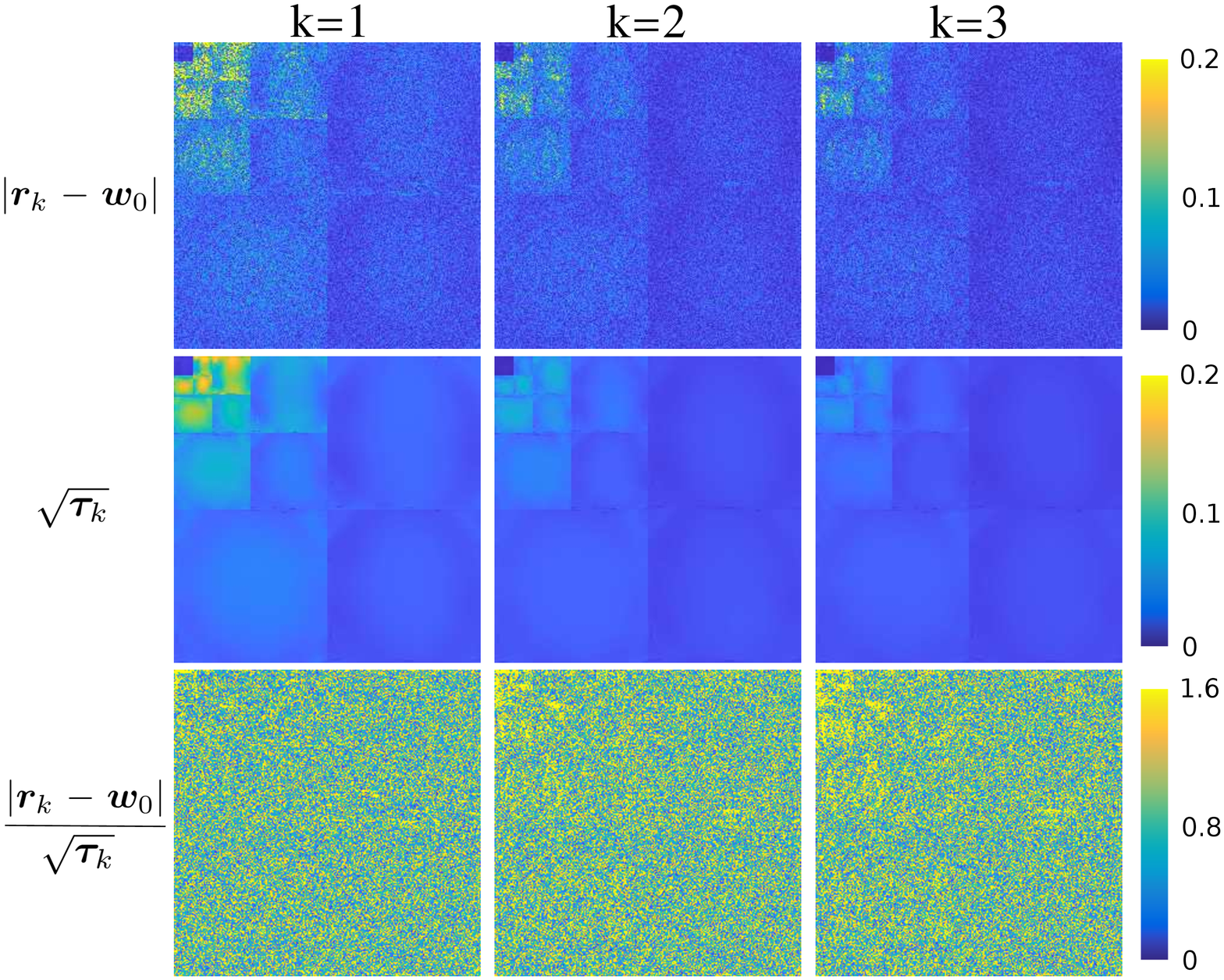}
\caption{\label{fig:tau_evo_k2} P-VDAMP's state evolution for the same brain dataset and sampling pattern as shown in Fig. \ref{fig:err_prop_example} at iterations $k=1, 2, 3$, showing that the Gaussian aliasing illustrated at $k=0$ approximately holds for subsequent iterations.}
\vspace{0.0cm}
\centering
\includegraphics[width=1 \columnwidth]{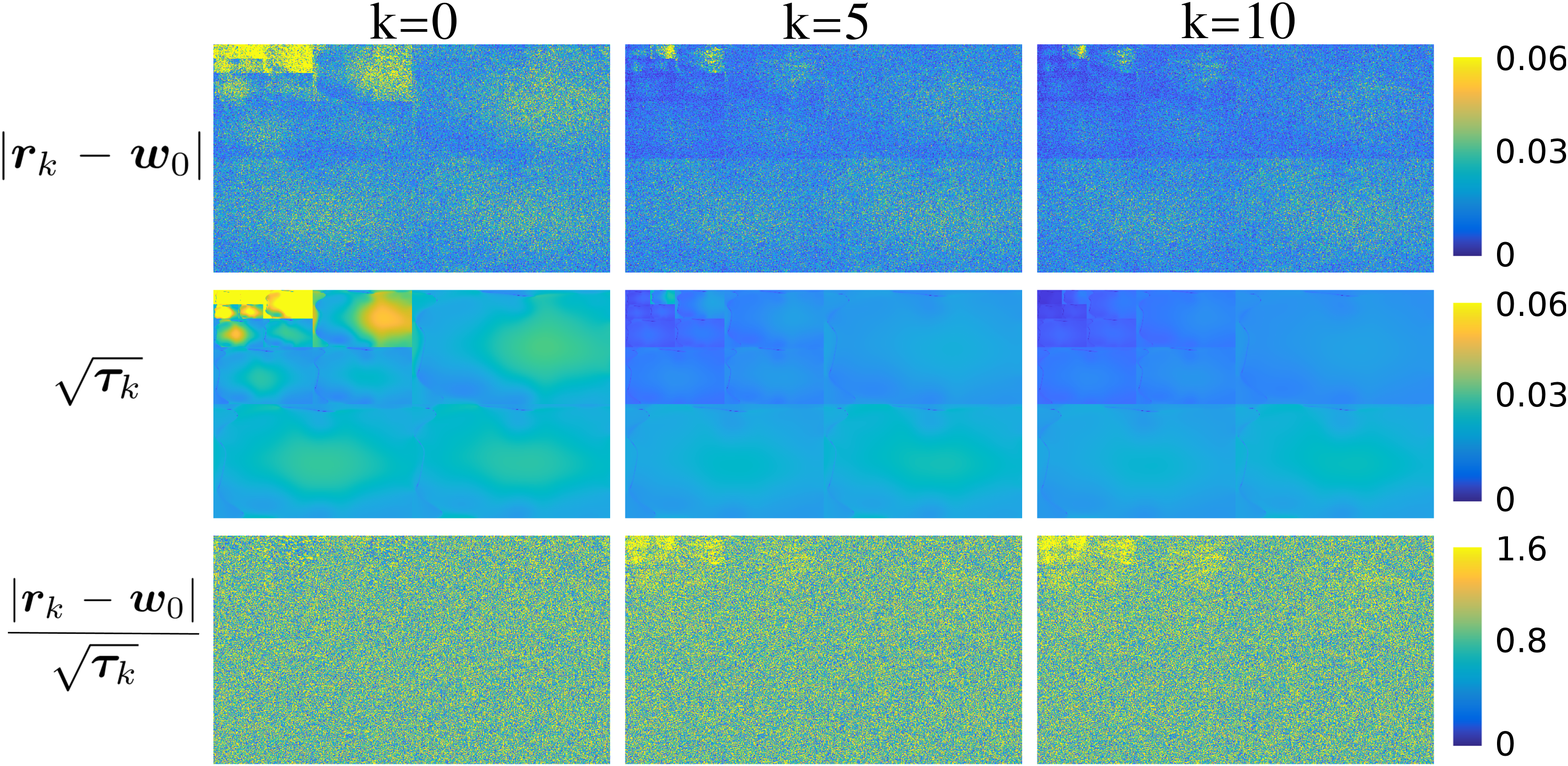}
\caption{\label{fig:tau_evo_knee} P-VDAMP's state evolution for a knee dataset sampled at $R=10$ at iterations $k=0, 5, 10$.  }
\vspace{-0.4cm}
\end{figure}

Fig. \ref{fig:brain_recon} shows one of the brain datasets at $R=5$. Optimal FISTA and both P-VDAMP estimates retain considerably more of the fine detail than competing methods. Shared FISTA's $\lambda$ is around $10$ larger than Optimal FISTA in this case, and compression artifacts are clearly visible. The competing parameter-free methods SURE-IT and A-FISTA peform better than Shared FISTA in this case, but still have substantial reconstruction artifacts.  The underlying anatomy is not  clearly visible in Unbiased P-VDAMP's error maps, with a similar contribution to the error at tissue boundaries and bulk tissue, which is a consequence of P-VDAMP's state evolution.

One of the knee datasets undersampled at $R = 10$ is shown in Fig. \ref{fig:knee_recon}. P-VDAMP outperforms both competing parameter-free methods. Shared FISTA's $\lambda$ is 6 times smaller than Optimal FISTA here, and performs reasonably well, but has a higher NMSE than P-VDAMP.  

Fig. \ref{fig:angio_recon} shows the reconstructions of one of the Angiograms at $R = 5$, where P-VDAMP outperforms Optimal FISTA. The green arrow in the magnified region highlights a narrow region of a blood vessel, which is visible for Unbiased P-VDAMP and partially visible for P-VDAMP but mostly lost in the competing methods.

The $\bm{\tau}_k$ update in line 5 of P-VDAMP, Algorithm \ref{alg:mult_coil_VDAMP}, is reasonably computationally expensive; P-VDAMP took around 4 times longer per iteration than FISTA in our MATLAB implementation. However, we found that P-VDAMP typically requires substantially fewer iterations, so required only a little longer than Optimal FISTA to reach convergence. This is illustrated in two NMSE vs time plots shown in Fig. \ref{fig:NMSE_vs_it_mc}.  

\subsection{P-VDAMP's state evolution \label{sec:mc_approx_se} } 

Fig. \ref{fig:tau_evo_k2} shows the wavelet-domain aliasing for one of the brain datasets at $k=1,2,3$ and one of the knee datasets at iterations $k=0,5,10$, both at $R=10$, which shows that P-VDAMP's parallel coloured state evolution \eqref{eqn:vdamp-se} approximately holds. There is a slight deviation from exactly unbiased aliasing apparent at larger $k$. This is due to P-VDAMP's damping factor $\rho$, which causes P-VDAMP's Onsager correction to be inexact. We found that using less damping caused state evolution to be more exact,  but that the algorithm converged to a poorer reconstruction. We therefore recommend using damping at the cost of a slight breakdown in state evolution.

\subsection{Robustness to measurement noise \label{sec:meas_noise_effect}}

To evaluate the effect of measurement noise on P-VDAMP's performance, we simulated Gaussian additive noise. We constructed an artificial measurement noise covariance matrix by generating a $N_c$-dimensional random vector $\bm{v}$ with complex entries with real and imaginary parts drawn from a uniform distribution between $-1$ and $1$, and scaled the vector to $|\bm{v}| = c \bm{1}_{N_c}$, where $c$ is a constant that determines the signal-to-noise ratio (SNR). For each measurement noise SNR we averaged over 10 random instances of $\bm{v}$.  We compared P-VDAMP with FISTA with a fixed sparse weighting, tuned to be optimal in the absence of simulated measurement noise. 

Fig. \ref{fig:noise_sens} shows the dependence of the reconstruction NMSE on measurement noise SNR for each algorithm for one of the knee datasets sampled at $R = 5, 10$. While FISTA's sparse weighting is fixed, P-VDAMP's sparse parameters adapt to the measurement noise on-the-fly and converges to a better solution for low SNR in these examples. Since the measurement noise is amplified by $\bm{P}^{-1}$ in P-VDAMP, the unbiased P-VDAMP output is more sensitive to measurement noise, and did not perform well for low SNR, indicating that it may not be appropriate for especially low SNR CS-MRI applications.

\section{Conclusions \label{sec:conc}}

P-VDAMP is an efficient, robust and principled tuning-free algorithm for CS-MRI that is competitive with Optimal FISTA and offers substantial robustness and reconstruction quality improvements over competing parameter-free methods.

State evolution essentially recasts multi-coil CS-MRI reconstruction as a sequence of Gaussian denoising tasks, where the Gaussian ``effective noise" has known distribution parameters. In other words, roughly speaking, how well one can reconstruct an image depends on how well one can denoise a Gaussian-corrupted version of it \cite{Metzler2016}. This suggests a number of potential advantages of P-VDAMP's state evolution beyond parameter tuning that could be explored in future work.  For instance, benefits for deep-learning based reconstruction \cite{NIPS2017_6774, borgerding2017amp, Metzler2020} and uncertainty quantification \cite{kitichotkul2020suremap} have already been explored in the context of single-coil VDAMP \cite{Metzler2020}. Initial work by the authors on Denoising-P-VDAMP (D-P-VDAMP) \cite{Milla2022, millard2021approximate} has shown that the deep Plug-and-Play framework \cite{Venkatakrishnan2013, Ahmad2020}  can improve P-VDAMP's performance.

We have found that P-VDAMP naively applied to Poisson disc sampled data also performs well, despite a violation of the independent sampling assumption. However, we have found that it typically diverges when applied to problems that contain smooth sampling trajectories, so cannot yet be applied to sampling schemes such as 3D radial \cite{Johnson2013} or 3D cones \cite{Gurney2006}, nor to purely 2D acquisitions. P-VDAMP could be extended to such sampling schemes by an appropriate modification of the aliasing model update in line 5 of Algorithm \ref{alg:mult_coil_VDAMP}.

Very recently an AMP-based algorithm was proposed \cite{shastri2022expectation} for single and multi-coil MRI acquisitions that uses similar principles to single-coil VDAMP but with an alternative procedure for data consistency at each iteration, which may have better fixed points. Extending the richer aliasing model presented in this paper to the multi-coil version of this alternative algorithm is a desirable avenue for future work.

\begin{appendix}
This appendix derives an approximate expression for the diagonal of the covariance of the aliasing of $\bm{r}_0 =\bm{\Psi }\hat{\bm{x}}_0$, where $\hat{\bm{x}}_0=\sum_c \bm{S}_c^H \bm{F}^H \bm{P}^{-1}\bm{y}_c$.  Let $\bm{y}_{0,c} =  \bm{F} \bm{S}_c\bm{x}_{0}$ and $\bm{w}_0 = \bm{\Psi}\bm{x}_0$ be the signal of interest's coil-weighted Fourier coefficients and wavelet coefficients respectively. Using $\bm{x}_0 =\sum_c \bm{S}_c^H \bm{F}^H \bm{y}_{0,c}$, the aliasing of $\bm{r}_0$ is
\begin{align*}
    |{\bm{r}}_0 - \bm{w}_0|^2 =&|\bm{\Psi}(\hat{\bm{x}}_0- \bm{x}_0)|^2  \\
    =& |\bm{\Psi}\sum_c \bm{S}_c^H \bm{F}^H[\bm{P}^{-1}\bm{y}_c - \bm{y}_{0,c}]|^2  \nonumber \\ 
    =& |\sum_c \hat{\bm{\Psi}}^c[\bm{P}^{-1}\bm{y}_c - \bm{y}_{0,c}]|^2,
\end{align*}
where $\hat{\bm{\Psi}}^c = \bm{\Psi}\bm{S}_c^H \bm{F}^H$. Since $\bm{y}_c =  \bm{M}_\Omega (\bm{y}_{0,c} + \bm{\varepsilon}_c)$,
\begin{align}
    & |\sum_c\hat{\bm{\Psi}}^c[\bm{P}^{-1}\bm{y}_c - \bm{y}_{0,c}]|^2 \nonumber \\ &= |\sum_c\hat{\bm{\Psi}}^c[(\bm{P}^{-1}\bm{M}_\Omega- \mathds{1})\bm{y}_{0,c} - \bm{P}^{-1}\bm{M}_\Omega \bm{\varepsilon}_c]|^2 \nonumber \\
    &= |\sum_c\hat{\bm{\Psi}}^c \bm{q}_{c}|^2 \label{eqn:Psixtild}
\end{align}
where we define $
\bm{q}_{c} = (\bm{P}^{-1}\bm{M}_\Omega- \mathds{1})\bm{y}_{0,c} - \bm{P}^{-1}\bm{M}_\Omega \bm{\varepsilon}_c.$ Computing the expectation of \eqref{eqn:Psixtild} with a method analogous to the single-coil approach in \cite{Milla2020} would involve computing and storing in memory the $ N \times N$ matrices $|\bm{\Psi}^c|^2$ for all $N_c$ coils, which is not practical for large $N$. We suggest approximating the $\hat{\bm{\Psi}}^c$ in terms of the unweighted spectra $\bm{\hat{\Psi}} = \bm{\Psi F}^H$, so that the spectral density of the $j$th row of $\hat{\bm{\Psi}}^c$ is expressed as 
\begin{equation}
\hat{\bm{\Psi}}^c_j \approx \xi_{c, j} \hat{\bm{\Psi}}_j, \label{eqn:xi_approx}
\end{equation}
where $\xi_{c, j} \in \mathds{C}$ is a constant.  Since $\hat{\bm{\Psi}}^c_j$ is the $j$th coil-weighted wavelet spectrum, \eqref{eqn:xi_approx} is interpretable as approximating the coil sensitivity as flat over the support of the $j$th wavelet filter, so that its spectrum is the unweighted spectrum multiplied by a filter- and coil-dependent scalar. 
We suggest using the $\xi_{c,j}$ that minimizes $\| \hat{\bm{\Psi}}_j^c - \xi \hat{\bm{\Psi}}_j \|^2_2$:
\begin{align*}
	\xi_{c,j} & = \argmin_{\xi} \| \hat{\bm{\Psi}}_j^c - \xi \hat{\bm{\Psi}}_j \|^2_2 \\
	& = \argmin_{\xi} \| \bm{\Psi}_j \bm{S}_c^H  \bm{F}^H - \xi \bm{\Psi}_j \bm{F}^H\|^2_2 \\
	& = \frac{\bm{\Psi}_j \bm{S}_c^H \bm{\Psi}_j^H}{\bm{\Psi}_j \bm{\Psi}_j^H} \\
	& = \braket{|\bm{\Psi}_j|^2, \mathrm{Diag}( \bm{S}_c^H)}, 
\end{align*}
where $\bm{\Psi}_j$ is the $j$th row of $\bm{\Psi}$, and we have used the orthonormality of the rows of $\bm{\Psi}_j$ to set $\bm{\Psi}_j \bm{\Psi}_j^H = 1$. 

\begin{figure}[t]
\centering
	\includegraphics[width =\columnwidth]{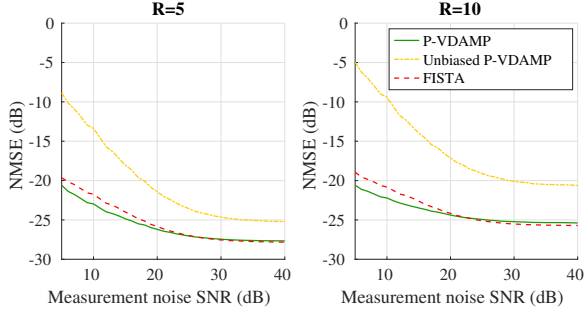}
      \caption{Reconstruction NMSE's dependence on simulated measurement noise SNR for one of the knee datasets. \label{fig:noise_sens}}
\end{figure}

\textbf{Claim 1.} 
	Let ${\xi}_{c,j} \in \mathds{C}$ be a constant that approximates $ \hat{\bm{\Psi}}_j^c \approx \xi_{c, j} \hat{\bm{\Psi}}_j $, for all $j$ rows of $\bm{\hat{\Psi}}^c$.  Then the expectation of \eqref{eqn:Psixtild} is
\begin{multline}
    \mathds{E} \{ |r_{0,j} - w_{0,j}|^2\} \\ \approx \bm{\xi}_{j}^H \left( \sum_i |\hat{\Psi}_{ji}|^2\left[ \left(\frac{1 - p_i}{p_i}\right) \bm{y}_{0,i} \bm{y}_{0,i}^H + \frac{1}{p_i}\bm{\Sigma}_{\varepsilon, i}^2  \right] \right) \bm{\xi}_{j}, \label{eqn:tau_expr_y0}
\end{multline}
where the $\bm{y}_{0,i} \in \mathds{C}^{N_c}$ and $\bm{\xi}_{0,j} \in \mathds{C}^{N_c}$ are comprised of the $i$th k-space location over all coils, having $c$th entry $[y_{0,i}]_c = [\bm{F S}_c \bm{x}_0]_i$ and $[\bm{\xi}_{0,j}]_c = \xi_{c, j}$,  and $\bm{\Sigma}_{\varepsilon, i}^2 = \mathds{E} \{ \bm{\varepsilon}_i \bm{\varepsilon}_i^H \}$ where $\bm{\varepsilon}_{i} \in \mathds{C}^{N_c}$ similarly has $c$th entry $\varepsilon_{i,c} = [\bm{\varepsilon}_c]_i$.
 
\begin{proof}
 Using \eqref{eqn:xi_approx}, the $j$th entry of \eqref{eqn:Psixtild} is approximately
\begin{align}
	\nonumber |{{r}}_{0,j} - {w}_{0,j}|^2 & \approx |\sum_{c}   \hat{\bm{\Psi}}_j \bm{q}_{c} \xi_{c, j} |^2 \\ & = |\sum_{c,i}   \hat{{\Psi}}_{ji} {q}_{c,i} \xi_{c, j} |^2. \label{eqn:todeltaj}
\end{align}
Using the result from Appendix C of \cite{Milla2020}, the expectation of \eqref{eqn:todeltaj} over the sampling mask and measurement noise is 
\begin{align}
	 \nonumber \mathds{E}\{|\sum_{c,i}   \hat{{\Psi}}_{ji} {q}_{c,i} \xi_{c, j} |^2\}  &= \sum_i |\hat{{\Psi}}_{ji}|^2 \mathds{E}\{| \sum_c  \xi_{c, j} q_{c,i} |^2\} \\
	&= \sum_i |\hat{{\Psi}}_{ji}|^2 \mathds{E}\{|  {\bm{\xi}_{j}^H} \bm{q}_{i} |^2\}, \label{eqn:err_prop_mid}
\end{align}
where $\bm{\xi}_{j} \in \mathds{C}^{N_c}$ and $\bm{q}_{i} \in \mathds{C}^{N_c}$.  Resolving the expectation,  
\begin{align}
   \nonumber \mathds{E} & \{|{\bm{\xi}_{j}^H} \bm{q}_{i} |^2\}  = \mathds{E}\{{\bm{\xi}_{j}^H} \bm{q}_{i} \bm{q}_{i}^H\bm{\xi}_{j} \} \\ =& \mathds{E}\left\{\bm{\xi}_{j}^H \left[ \left(\frac{m_i}{p_i} - 1 \right)^2 \bm{y}_{0,i} \bm{y}_{0,i}^H   +\left(\frac{m_i}{p_i}\right)^2\bm{\varepsilon}_i \bm{\varepsilon}_i^H \right]\bm{\xi}_{j} \right\}, \label{eqn:bob}
\end{align}
where $\mathds{E}\{ \bm{\varepsilon_i}\} = \bm{0}_{N_c}$ has been applied to evaluate the cross-terms to zero. For Bernoulli sampling,
\begin{equation*}
\mathds{E}\left\{\left(\frac{m_i}{p_i} - 1 \right)^2\right \} = \frac{1-p_i}{p_i}, 
\end{equation*}
Therefore, defining $\mathds{E} \{ \bm{\varepsilon}_i \bm{\varepsilon}_i^H \} = \bm{\Sigma}_{\varepsilon, i}^2$, \eqref{eqn:bob} is
\begin{align}
  \mathds{E} \{|\bm{\xi}_{j}^H \bm{q}_{i} |^2\}  = \bm{\xi}_{j}^H\left[ \left(\frac{1 - p_i}{p_i}\right) \bm{y}_{0,i} \bm{y}_{0,i}^H + \frac{1}{p_i}\bm{\Sigma}_{\varepsilon, i}^2 \right] \bm{\xi}_{j}. \label{eqn:tau_expry0}
\end{align}
Substituting \eqref{eqn:tau_expry0} in to \eqref{eqn:err_prop_mid} gives \eqref{eqn:tau_expr_y0} as required.
\end{proof} 

Eqn. \eqref{eqn:tau_expr_y0} requires knowledge of $\bm{y}_{0}$ so is of limited practical use. However, we can approximate $\mathds{E}_{\Omega,\varepsilon} \{ |r_{0,j} - w_{0,j}|^2\}$ without the ground truth with
\begin{equation}
	 \tau_j = \bm{\xi}_j^H \left( \sum_i |\hat{\Psi}_{ji}|^2\left[ \left(\frac{1 - p_i}{p_i}\right) \bm{y}_{i} \bm{y}_{i}^H + \bm{\Sigma}_{\varepsilon, i}^2 \right]\frac{m_i}{p_i} \right) {\bm{\xi}_{j}}. 
\label{eqn:tau_full}
\end{equation}

 \textbf{Claim 2.}  \label{clm:imp_samp}
	The estimate of the multi-coil wavelet-domain error	stated in \eqref{eqn:tau_full} is unbiased:
	\begin{equation}
    \mathds{E} \{\tau_j\} = \bm{\xi}_{j}^H \left( \sum_i |\hat{\Psi}_{ji}|^2\left[ \left(\frac{1 - p_i}{p_i}\right) \bm{y}_{0,i} \bm{y}_{0,i}^H + \frac{1}{p_i}\bm{\Sigma}_{\varepsilon, i}^2 \right] \right) {\bm{\xi}_{j}}. \label{eqn:tau_expr}
	\end{equation}
\begin{proof} Since $\bm{y}_{i} = m_i (\bm{y}_{0,i} + \bm{\varepsilon}_i)$, and $\mathds{E}\{ \bm{\varepsilon}_i\} = \bm{0}_{N_c}$,
\begin{align}
 & \mathds{E} \left \{ \left[ \left(\frac{1 - p_i}{p_i}\right) \bm{y}_{i} \bm{y}_{i}^H + \bm{\Sigma}_{\varepsilon, i}^2  \right]\frac{m_i}{p_i} \right \} \nonumber \\ & =    \mathds{E} \left \{ \left[ m_i \left(\frac{1 - p_i}{p_i}\right) (\bm{y}_{0,i}  + \bm{\varepsilon}_i)(\bm{y}_{0,i}  + \bm{\varepsilon}_i)^H  + \bm{\Sigma}_{\varepsilon, i}^2 \right]\frac{m_i}{p_i} \right \} \nonumber \\
 & =  \left(\frac{1 - p_i}{p_i}\right) \bm{y}_{0,i} \bm{y}_{0,i}^H + \frac{1}{p_i} \bm{\Sigma}_{\varepsilon, i}^2 \label{eqn:mc_imp_sampling}
\end{align}
Substituting \eqref{eqn:mc_imp_sampling} in to the expectation of \eqref{eqn:tau_full} gives the right-hand-side of \eqref{eqn:tau_expr}.
\end{proof}
\end{appendix}

\bibliographystyle{ieeetr}
\bibliography{library_manu}

\end{document}